\documentclass[12pt]{article}

\usepackage{epsf,graphicx, color}
\usepackage{tikz}
\usetikzlibrary{arrows,automata}
\usepackage{cite, amsmath, amssymb}
\usepackage{color}
\newtheorem{theorem}{\bf Theorem}[section]
\newtheorem{corollary}[theorem]{\bf Corollary}

\newtheorem{proposition}[theorem]{\bf Proposition}

\newcommand{\proof}{\noindent{\bf Proof.\ }}
\newcommand{\qed}{\hfill $\square$ \bigskip}
\newcommand{\ecc}{{\rm ecc}}

\newcommand{\diam}{{\rm diam}}
\newcommand{\Tr}{{\rm Tr}}
\newcommand{\cp}{\,\square\,}

\begin{document}

\title{Generalized stepwise transmission irregular graphs}

\author{Yaser Alizadeh$^{a}$, Sandi Klav\v zar$^{b,c,d}$, Zohre Molaee$^a$}
\date{}
\maketitle
% \vspace{-0.8 cm}
\begin{center}

$^{a}$ Department of Mathematics, Hakim Sabzevari University, Sabzevar, Iran \\
{\tt  y.alizadeh@hsu.ac.ir}  \\
{\tt  zohremolaeegm@gmail.com }  \\
$^b$ Faculty of Mathematics and Physics, University of Ljubljana, Slovenia\\
{\tt sandi.klavzar@fmf.uni-lj.si} \\
\medskip

$^c$ Institute of Mathematics, Physics and Mechanics, Ljubljana, Slovenia \\
\medskip

$^d$ Faculty of Natural Sciences and Mathematics, University of Maribor, Slovenia 
\end{center}

\begin{abstract} 
The transmission ${\rm Tr}_G(u)$ of a vertex $u$ of a connected graph $G$ is the sum of distances from $u$ to all other vertices. $G$ is a stepwise transmission irregular (STI) graph if $|{\rm Tr}_G(u) - {\rm Tr}_G(v)|= 1$ holds for any edge $uv\in E(G)$.  In this paper, generalized STI graphs are introduced as the graphs $G$ such that for some $k\ge 1$ we have $|{\rm Tr}_G(u) - {\rm Tr}_G(v)|= k$ for any edge $uv$ of $G$. It is proved that generalized STI graphs are bipartite and that as soon as the minimum degree is at least $2$, they are 2-edge connected. Among the trees, the only generalized STI graphs are stars. The diameter of STI graphs is bounded and extremal cases discussed. The Cartesian product operation is used to obtain highly connected generalized STI graphs. Several families of generalized STI graphs are constructed. 
\end{abstract}


\medskip\noindent
{\bf AMS Subj.\ Class:} 05C12

%%%%%%%%%%%%%%%%%%%%%%%%%%%%%%%%%%%%%%%
\section{Introduction}
\label{sec:intro}

The shortest-path distance $d_G(u,v)$ between vertices $u$ and $v$ of a graph $G$ is the minumum number of edges on a $u,v$-path. The {\em transmission} $\Tr_G(u)$ of a vertex $u$ is the sum of distances between $u$ and all the other vertices in $G$. Two early papers in which the transmission was considered are~\cite{P-1984, S-1991}, where the interest was on maximal transmission in several classes of graphs and on the behaviour of the transmission under removing a vertex. Transmission plays an important role in the investigation of distance-based graph invariants such as the the Wiener index~\cite{knor-2016} and the Mostar index~\cite{ali-2021}. In particular, several measures on transmission irregularity were posed in~\cite{ali-2021}. The fact that the transmission is a fundamental concept in metric graph theory and wider is demonstrated by the fact that it is also known by other names such as the total distance of a vertex~\cite{cava-2019, kla-2013} and the status of a vertex~\cite{abiad-2021, QZ2020}. 

Interesting graph families have recently been defined based on the transmission. {\em Transmission irregular graphs} are the graphs in which any two different vertices $G$ have different transmissions. Although (or perhaps because) almost no graph is transmission irregular~\cite{AK-2018}, the search for such graphs has become of interest to several groups of researchers, some of the selected papers on this topic are~\cite{AS-2020, BD-2021, D-2019, XK-2021}. If we further require that the vertex transmissions of a graph form a sequence of consecutive integers, then we speak of an {\em interval transmission irregular graph}~\cite{al-yakoob-2022a}. 

{\em Stepwise transmission irregular graphs}, {\em STI graphs} for short, are the graphs in which for every edge the transmissions of its endpoints differ by $1$. STI graphs were introduced in~\cite{DS-2020}. The research was continued in~\cite{al-yakoob-2022b} where a conjecture from~\cite{DS-2020} was confirmed that all graphs from a certain family are STI. Moreover, a computational support was provided for another conjecture from~\cite{DS-2020} asserting that each STI graph has girth $4$. In general, however, the conjecture  remains open. STI graphs which are extremal with respect to different metric invariants such as the diameter, the Wiener index, and the eccentricity index, were characterized in~\cite{AK-2023}. 
 
\begin{figure}[ht!]
\begin{center}
\begin{tikzpicture}[scale=1,style=thick]
\tikzstyle{every node}=[draw=none,fill=none]
\def\vr{2.5pt} 

\begin{scope}[yshift = 0cm, xshift = 0cm]
\path (5,4) coordinate (u1);
\path (6.5,4) coordinate (u2);
\path (8,4) coordinate (u3);
\path (8,3) coordinate (u4);
\path (8,2) coordinate (u5);
\path (6.5,2) coordinate (u6);
\path (5,2) coordinate (u7);
\path (5,3) coordinate (u8);
\path (6.5,3.5) coordinate (u9);
\path (6,3) coordinate (u10);
\path (11,4) coordinate (w1);
\path (12.5,4) coordinate (w2);
\path (14,4) coordinate (w3);
\path (14,3) coordinate (w4);
\path (14,2) coordinate (w5);
\path (12.5,2) coordinate (w6);
\path (11,2) coordinate (w7);
\path (11,3) coordinate (w8);
\path (12.5,3) coordinate (w9);

%% edges %%
\draw (u1) -- (u2) -- (u3) -- (u4)  -- (u5) -- (u6) -- (u7) -- (u8) --(u1);
\draw (u2) --(u6);
\draw (u4) -- (u8);
\draw (w1) -- (w2) -- (w3) -- (w4)  -- (w5) -- (w6) -- (w7) -- (w8) --(w1);
\draw (w2) -- (w9) -- (w6);
\draw (w4) -- (w9) -- (w8);
%% vertices %%%
\draw (u1)  [fill=white] circle (\vr);
\draw (u2)  [fill=white] circle (\vr);
\draw (u3)  [fill=white] circle (\vr);
\draw (u4)  [fill=white] circle (\vr);
\draw (u5)  [fill=white] circle (\vr);
\draw (u6)  [fill=white] circle (\vr);
\draw (u7)  [fill=white] circle (\vr);
\draw (u8)  [fill=white] circle (\vr);
\draw (u9)  [fill=white] circle (\vr);
\draw (u10)  [fill=white] circle (\vr);
\draw (w1)  [fill=white] circle (\vr);
\draw (w2)  [fill=white] circle (\vr);
\draw (w3)  [fill=white] circle (\vr);
\draw (w4)  [fill=white] circle (\vr);
\draw (w5)  [fill=white] circle (\vr);
\draw (w6)  [fill=white] circle (\vr);
\draw (w7)  [fill=white] circle (\vr);
\draw (w8)  [fill=white] circle (\vr);
\draw (w9)  [fill=white] circle (\vr);
%% text %%
\draw[above] (u1) node {$20$}; 
\draw[above] (u2) node {$18$}; 
\draw[above] (u3) node {$20$}; 
\draw[right] (u4) node {$18$}; 
\draw[below] (u5) node {$20$}; 
\draw[below] (u6) node {$18$}; 
\draw[below] (u7) node {$20$}; 
\draw[left] (u8) node {$18$}; 
\draw[right] (u9) node {$20$}; 
\draw[below] (u10) node {$20$}; 
\draw[above] (w1) node {$18$}; 
\draw[above] (w2) node {$15$}; 
\draw[above] (w3) node {$18$}; 
\draw[right] (w4) node {$15$}; 
\draw[below] (w5) node {$18$}; 
\draw[below] (w6) node {$15$}; 
\draw[below] (w7) node {$18$}; 
\draw[left] (w8) node {$15$}; 
\draw[right] (w9) ++(0.0,-0.3) node {$12$}; 
\end{scope}

\end{tikzpicture}
\end{center}
\caption{A $2$-STI graph (left) and a $3$-STI graph (right)}
\label{fig1}
\end{figure}
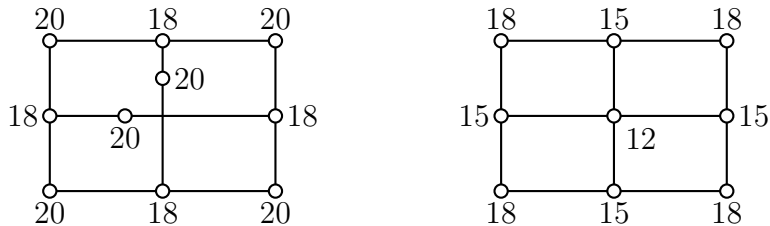

In this paper we extend STI graphs to generalized STI graphs as follows. If $k$ is an arbitrary positive integer, then we say that a graph is a {\em $k$-STI graph} if for every edge the transmissions of its endpoints differ by $k$. If $G$ is a $k$-STI graph for some $k$, then we say that $G$ is a {\em generalized STI graph}. A $2$-STI graph and a $3$-STI graph are shown in Fig.~\ref{fig1}, where next to each vertex its transmission is written. 
 
We proceed as follows. After this paragraph, some definition and a useful result needed are stated. In the next section several properties of generalized STI are established. Among other results we prove that generalized STI graphs are bipartite, and that as soon a generalized STI graphs has $\delta(G)\ge 2$, it is 2-edge connected. We also bound the diameter of generalized STI graphs and deduce that among trees the only generalized STI graphs are stars. In Section~\ref{sec:graph operation} we prove that if $G$ and $H$ have order $n$, then $G \cp H$ is an $(nk)$-STI graph if and only if $G$ and $H$ are $k$-STI graphs. Finally, In Section~\ref{sec:constructions}, several families of generalized STI graphs are constructed.

The graphs considered here are simple and connected. Let $G =(V(G),E(G))$ be a graph. Then the order of $G$, the degree of a vertex $v$, and the minimum degree of $G$ are denoted by $n(G)$, $\deg(v)$, and $\delta(G)$ respectively. 

A subgraph $H$ of a graph $G$ is {\em isometric} if $d_H(u,v) = d_G(u,v)$ holds for each pair of vertices $u,v\in V(H)$. The {\em eccentricity} $\ecc_G(u)$ of a vertex $u\in V(G)$ is the maximum distance from $u$ to other vertices of $G$. The {\em diameter} $\diam(G)$ of $G$ is the maximum eccentricity of vertices in $G$. For an edge $e=uv\in E(G)$, let  $N_u(e)$ denote the set of vertices of $G$ that are closer to $u$ than to $v$. Similarly $N_v(e)$ is defined. Let further $n_u(e) = |N_u(e)|$ and $n_v(e) = |N_v(e)|$. If the edge $e$ lies in different graphs and we consider it as an edge of a graph $G$, then we will specify this notation to $n_u(e|G)$ and $n_v(e|G)$. As already said, $\Tr_G(u) = \sum_{x\in V(G)} d_G(u,x)$. For an edge $e = uv$ we will further use the notation $I(e) = |\Tr_G(u) - \Tr_G(v)|$ and call $I(e)$ the {\em transmission imbalance} of the edge $e$.  

The following useful result goes back to~\cite{Entringer}. 

\begin{theorem}
\label{thm:Entringer}
If $G$ is a bipartite graph and $e=uv\in E(G)$, then $\Tr_G(u)  - \Tr_G(v)= n_v(e) - n_u(e)$. 
\end{theorem}

%%%%%%%%%%%%%%%%%%%%%%%%%%%%%%%%%%%%%%%%%%%%%%%%%%%%%%
\section{Properties of $k$-STI graphs}
\label{sec:properties}
 
Consider the complete bipartite graph $K_{p,q}$, where $p\ge q\ge 1$. Then the transmissions of the vertices from the two bipartition sets are $q + 2p$ and $p + 2q$, respectively. Hence, if $p>q$, then $K_{p,q}$ is a $(p-q)$-STI graph. For a given $k\ge 1$, the graphs $K_{n+k,n}$, $n\ge 1$, thus form an infinite family of $k$-STI graphs. 

Our first result collects the structural properties of $k$-STI graphs. To this end, recall that vertices $u$ and $v$ are {\em twins} if for every $w\in V(G)\setminus \{u,v\}$ we have $uw\in E(G)$ if and only if $vw\in E(G)$. Note that in our definition of twin vertices we put no condition whether $u$ and $v$ are adjacent.  We say that $G$ is {\em twin-free}, if it contains no twins.    

\begin{theorem}\label{thm:base}
If $k\ge 1$ and $G$ is a $k$-STI graph, then the following hold. 
\begin{enumerate}
\item[(i)] $G$ is a bipartite, twin-free graph.  
\item[(ii)] $n(G) \equiv k \pmod{2}$.
\item[(iii)] If $\delta(G)\ge 2$, then $G$ is $2$-edge connected.
\item[(iv)] If $k\in [2]$, then $G$ is either $2$-connected or $G\in \{P_3, K_{1,3}\}$.
\end{enumerate}
\end{theorem}

\proof
(i) Suppose on the contrary that $G$ contains an odd cycle $v_0 v_1 \ldots v_{2p}$. As $G$ is a $k$-STI graph, we have  $\Tr_G(v_0) \equiv \Tr_G(v_{2}) \pmod{2k}$, $\ldots$,     $\Tr_G(v_{2p-2}) \equiv \Tr_G(v_{2p}) \pmod{2k}$, so that also $\Tr_G(v_0) \equiv \Tr_G(v_{2p}) \pmod{2k}$. Since $v_0v_{2p}\in E(G)$, this is not possible. Hence $G$ is bipartite. It is straightforward to see that if $u$ and $v$ are twins, then $\Tr_G(u) = \Tr_G(v)$, thus $G$ must also be twin-free.  

(ii) Let $e = uv\in E(G)$, and assume without loss of generality that $\Tr_G(u) = \Tr_G(v) + k$.  As just proved in (i), $G$ is bipartite, hence Theorem~\ref{thm:Entringer} yields that $\Tr_G(u) - \Tr_G(v) = k = n_u(e) - n_v(e)$. Using once more that $G$ is bipartite we also have $n_u(e) + n_v(e) = n(G)$. It follows that $n(G) + k = 2n_u(e)$ which means that  $n(G) \equiv k \pmod{2}$.  

(iii) Suppose on the contrary that $e=uv$ is a bridge of $G$. Let $G-e = G_1 \cup G_2$, where $u \in V(G_1)$ and $v \in V(G_2)$. We may assume without loss of generality that $n_1 = n(G_1) \ge n_2 = n(G_2)$. By (i) we know that $G$ is  bipartite, hence using Theorem~\ref{thm:Entringer} we have 
$$k = I(uv) = n_u(e) - n_v(e) = n_1 - n_2\,.$$ 
Since $\delta(G)\ge 2$ and so also $\deg_G(v) \ge 2$, there exists a vertex $w\in V(G_2)$ such that $e' = vw\in E(G_2)$.  Then  
$$I(vw) = n_v(e') - n_w(e') \ge (n_1 +1) -(n_2 -1) = k+2\,,$$ 
a contradiction.

(iv) In~\cite[Proposition 5]{DS-2020} is was proved that all 1-STI graphs but the path $P_3$ are $2$-connected. Hence it remains to consider the case $k=2$. It is straightforward to check that $K_{1,3}$ is the only $2$-STI graph of order at most $4$. Assume in the rest that $G$ is a $2$-STI graph with $n = n(G) \ge 5$. If $G$ has a pendant edge $uv$ then $2 = |\Tr_G(u) - \Tr_G(v)|= n-2$ which is not possible since $n\ge 5$. Thus $\delta(G) \ge 2$. Suppose on the contrary that $G$  contains a cut vertex $v$. Let $G_1$ be a component of $G-v$ with the minimum order $n_1 = n(G_1)$. Notice that $n \ge 2n_1 +1$. As $\delta(G)\ge 2$, we have $n_1\ge 2$, and let  $w,z \in V(G_1)$ such that $vw, zw \in E(G)$. By (i)  $G$ is bipartite and hence $z$ is not adjacent to $v$. Thus 
$$2= \Tr_G(z) - \Tr_G(w) \ge (n-n_1+1) - (n_1-1) = n-2n_1 + 2 \ge 3\,.$$ 
This contradiction completes the argument. 
\qed

The assertions (i) and (ii) of Proposition~\ref{thm:base} respectively extend \cite[Proposition 1]{DS-2020} and \cite[Proposition 2]{DS-2020} which assert that STI graphs are bipartite graphs of odd order. The two examples of Fig.~\ref{fig1} imply that the assertion (iii) cannot be extended in general to $\ell$-edge connectedness for $\ell \ge 3$, while the graphs $K_{k+1,1}$ demonstrate that the assumption $\delta\ge 2$ cannot be avoided. Finally, (iv) also does not extend to $k\ge 3$. For instance, consider the graph obtained from two disjoint $4$-cycles by identifying a vertex from each of them. This graph is a $3$-STI graph with a cut vertex. For further such examples see  Proposition~\ref{prop:not2-con}.  

We next show: 

\begin{theorem}\label{max-order}
If $G$ is a $k$-STI graph, $k\ge 1$, then $k \le n(G) -2$. Moreover, the equality holds if and only if $G \cong K_{1,k+1}$. 
\end{theorem}

\proof
Let $e=uv$ be an edge of $G$. Since $G$ is bipartite by Theorem~\ref{thm:base}(i), we may assume without loss of generality (having Theorem~\ref{thm:Entringer} in mind) that $n_u(e) > n_v(e)$, so that $n_u(e) - n_v(e) = k$ holds. Since $n_u(e) +n_v(e) = n(G)$,  we infer that $n(G) = k + 2n_v(e) \ge k+2$ with equality holding if and only if $n_v(e) = 1$. This implies that $v$ is a pendant vertex. If $w$ is another vertex adjacent to $u$, then  $|\Tr_G(u) - \Tr_G(w)|= k = n(G) -2$. Thus $w$ is also a leaf. We conclude that all vertices adjacent to $u$ are pendant which in turn implies that $G\cong K_{1,k+1}$.
\qed

\begin{corollary}\label{tree}
If $k\ge 1$, then a tree $T$ is a $k$-STI graph if and only if $T \cong K_{1,k+1}$. 
\end{corollary}
 
\proof
If $e= uv\in E(T)$, where $v$ is a leaf of $T$, then $I(uv) =n_u(e) -n_v(e) = n(T)-2$. The result now follows from Theorem~\ref{max-order}. 
\qed

To conclude the section we bound the diameter of $k$-STI graphs as follows. 

\begin{theorem}\label{diameter}
If $G$ is an $k$-STI graph of order $n\ge 5$, then
$$ 2\le \diam(G) \le \frac{n+k}{2}-1\,.$$
Moreover, the left equality holds if and only if $G \cong K_{\frac{n+k}{2},\frac{n-k}{2}}$.
\end{theorem} 

\proof
Since $G$ is bipartite and $nG)\ge 5$, we have $\diam(G)\ge 2$. Moreover,  $\diam(G) = 2$ if and only if $G$ is a complete bipartite graph $K_{p,q}$, $p > q$. But then $p + q = n$ and $p-q = k$ which yields that $G \cong K_{\frac{n+k}{2},\frac{n-k}{2}}$. This proves the lower bound and the equality case. 

To prove the upper bound, we claim that $G$ contains an edge $g = xy$ such that $n_{x}(g) \ge \diam(G) + 1$.  Let $P$ be a diametral path in $G$ with $v$ and $w$ its endpoints. Let $v'$ be the neighbor of $v$ on $P$ and let $w'$ be the neighbor of $w$ on $P$. We may assume that $v'\ne w'$, for otherwise $\diam(G) = 2$ and the upper bound clearly holds. Let $e = vv'$ and $f=ww'$. If $\deg(v) = 1$, then $n_{v'}(e) \ge \diam(G) + 1$ because $n\ge 5$ and therefore $G$ is not a path. Similarly, $n_{w'}(e) \ge \diam(G) + 1$ holds if $\deg(w) = 1$. Hence assume that $\deg(v)\ge 2$ and $\deg(w)\ge 2$. Let $w''$ be a neighbor of $w$, $w''\ne w'$. Note that all the vertices of $P$ but $v$ lie in $N_{v'}(e)$. Hence, if also $w''\in N_{v'}(e)$, then  $n_{v'}(e) \ge \diam(G) + 1$. Assume hence that $w''\notin N_{v'}(e)$. Then $d_G(v',w'') = \diam(G)$ and $d_G(v,w'') = \diam(G) - 1$. Let $Q$ be a shortest $v,w''$-path. If $v'$ lies on $Q$, then $Q$ contains another vertex which is not on $P$ and lies in $N_{v'}(e)$, hence $n_{w'}(e) \ge \diam(G) + 1$ as required. Assume next that the neighbor $v''$ of $v$ on $Q$ is different from $v'$. Moreover, all the vertices on $Q$ lie in $N_{v}(e)$.  Now, since $G$ is not a cycle, there exists a vertex $z\notin V(P)\cup V(Q)$. As $G$ is bipartite, either $z\in N_{v}(e)$ or $z\in N_{v}(e)$. In either case, the existence of a required edge is proved. 

Assume now without loss of generality that $n_{v}(e) \ge \diam(G) + 1$. Since $G$ is an $k$-STI graph, we have $|n_v(e) - n_{v'}(e)| = k$ and $n_v + n_{v'} = n$. Thus $\diam(G) \le n_v (e) -1 \le \frac{n+k}{2}-1$. 
\qed

In~\cite[Lemma 3.1]{AK-2023} an infinite family of graphs was constructed for which the upper bound in Theorem~\ref{diameter} is attained. It would be of interest to construct such families for each $k\ge 2$ (or prove they do not exist). A sporadic example for $k=3$ is the graph obtained from two $4$-cycles by identifying a vertex from one by a vertex from the other. 

%%%%%%%%%%%%%%%%%%%%%%%%%%%%%%%%%%%%%%%%%%%%%%%%%%%%%%%%%%%
\section{Graph operations and $k$-STI graphs}
\label{sec:graph operation}

Using Theorem~\ref{thm:base} one can show that many local or global graph operations do not preserve the property of being generalized STI. For instance, by Theorem~\ref{thm:base}(i), the line graph $L(G)$ of a generalized STI graph $G$ is not such except $L(P_4) = P_3$. Similarly, by checking the small cases and by applying Theorem~\ref{thm:base}(i), the complement of a generalized STI graph $G$ is never a generalized STI graph. 

%%%%%%%%%%%%%%%%%%%%%%%
\iffalse
Another such situation is the following. 

\begin{proposition}
\label{prop:not-preserving}
If $G$ is a $k$-STI graph, then and if $e\in E(G)$ is not a bridge, then for any $\ell$, the graph $G-e$ is not a $\ell$-STI graph.
\end{proposition}

\proof
Let $e=vw$  and set $G' = G-e$. As $e=v$ is not a bridge, $e$ lies in at least one cycle. Let $C: v= v_1, v_2, \cdots, v_{2t} = w, v_1$ be a cycle of minimum length containing $e$. Then $C$ is an isometric subgraph of $G$. Then $\Tr_{G'}(v_t) = \Tr_{G}(v_t)$ and $\Tr_{G'}(v_{t+1}) = \Tr_{G}(v_{t+1})$. Hence, if $G'$ is a $\ell$-STI graph, then $\ell = k$. Assume without loss of generality that $\Tr_{G'}(v_t) = \Tr_{G'}(v_{t+1}) + k$. But then since $\Tr_{G'}(v_{t+2}) > 
$

In particular, if  $x=v_{\frac{t}{2}+1}$, then $ d_{G'}(x,v)=d_G(x,v)$ and $d_{G'}(x,w)=d_G(x,w)$. This yields $\Tr_{G'}(x) = \Tr_{G}(x)$. For vertex $z=v_{\frac{t}{2}+2}$ similar result is obtained.
But for vertex $y=v_{\frac{t}{2}}$ we get $d_{G'}(y,v) \ge d_G(y,v)+2$ and $d_G(y,w) = d_{G'}(y,w)$. Moreover
$\Tr_{G'}(y) \ge \Tr_{G}(y)+2$. Thus $|\Tr_{G'}(y)-\Tr_{G'}(x)| \ge k + 2$, hence $G-e$ is not a $k$-STI graph. 
\qed
\fi
%%%%%%%%%%%%%%%%%%%%%%%

In the previous section we observed that the graphs $K_{n+k,n}$, $n\ge 1$, are $k$-STI graphs. To obtained more involved highly connected generalized STI graphs, the Cartesian product can be used. Recall that the {\em Cartesian product} of two graphs $G$ and $H$, denoted $G\cp H$, is the graph with vertex set $V(G\cp H) = V(G) \times V(H)$ and vertices $(u, v)$ and $(x, y)$ are adjacent in $G\cp H$ if either $u=x$ and $vy \in E(H)$ or $v=y$ and $ux \in E(G)$, see~\cite{HIK-2011} for more information on this graph operation. 

\begin{theorem}\label{thm:cartesian}
Let $G$ and $H$ be graphs with $n(G) = n(H) = n$ and let $k\ge 1$. Then $G \cp H$ is an $(nk)$-STI graph if and only if $G$ and $H$ are $k$-STI graphs.
\end{theorem}

\proof 
From~\cite{AA-2014} we recall that if $x\in V(G)$ and $u\in V(H)$, then 
\begin{equation}
\label{eq:cp}
\Tr_{G \cp H}((x,u)) = \Tr_G (x) n(H) + \Tr_H (u) n(G)\,.
\end{equation}
Assume first that $G \cp H$ is an $nk$-STI graph. If $xy\in E(G)$ and $u\in V(H)$, then $(x,u)(y,u)\in E(G\cp H)$ and thus~\eqref{eq:cp} yields 
\begin{align*}
nk & = |\Tr_{G \cp H}((x,u)) - \Tr_{G \cp H}((x,v))| \\
& =  |\Tr_G (x) n + \Tr_H (u) n - \Tr_G (y) n -\Tr_H (u) n| \\
& = n |\Tr_G (x) - \Tr_G (y)|\,, 
\end{align*}
and hence $|\Tr_G (x) - \Tr_G (y)| = k$. IT follows that $G$ is a $k$-STI graph. Analogously we see that $H$ is $k$-STI graph.

Conversely, assume that $G$ and $H$ are $k$-STI graphs. If $(x,u)$ and $(x,v)$ are adjacent vertices in $G \cp H$, then applying~\eqref{eq:cp} once more we have
\begin{align*}
|\Tr_{G \cp H}(x,u)-\Tr_{G \cp H}(x,v)| & = |\Tr_G (x)n + \Tr_H(u)n - \Tr_G(x)n - \Tr_H(v)n\\
 & = n |\Tr_H(u) - \Tr_H(v)| = nk\,.
\end{align*}
Analogously we get the same conclusion for the edges of $G\cp H$ whose endvertices differ in the first coordinate. We conclude that $G\cp H$ is an $nk$-STI graph.
\qed

If $G$ and $H$ are graphs on at least two vertices, then the following formula applies to the connectivity the connectivity $\kappa (G\cp H)$ of $G\cp H$:  
$$\kappa (G\cp H) = \min\{ \kappa(G) n(H), \kappa(H) n(G), \delta (G) + \delta (H)\}\,.$$
The formula was announced in 1978 in~\cite{LI-1978}. However, neither its proof was provided nor did it appear afterwards. After several partial results, the formula was proved in 2008 by \v{S}pacapan in~\cite{SP-2008}. An appealing consequence of the formula is that $\kappa (G\cp H) \geq \kappa(G) + \kappa(H)$ holds for any connected graphs $G$ and $H$, cf.~\cite[Exercise 25.4]{HIK-2011}. Another consequence of the formula is that if $G$ is a connected graph of order at least $2$ and $G^{\cp, n}$ denotes the Cartesian product of $n$ copies of $G$, then for any $n\ge 2$ we have $\kappa (G^{\cp, n}) = \delta (G^{\cp, n}) = n\, \delta (G)$~\cite{KLSP-08}. These two consequences together with Theorem~\ref{thm:cartesian} guarantee the existence of numerous highly connected generalized STI graphs. 

%%%%%%%%%%%%%%%%%%%%%%%%%%%%%%%%%%%%%%%%%%%%%%%%%%
\section{Some families of generalized STI graphs}
\label{sec:constructions}

In this section, some families of generalized STI graphs are constructed. 

\begin{figure}[ht!]
\begin{center}
\begin{tikzpicture}[scale=0.6,style=thick]
\tikzstyle{every node}=[draw=none,fill=none]
\def\vr{2.5pt} 

\begin{scope}[yshift = 0cm, xshift = 0cm]
\path (6,7) coordinate (v1);
\path (6.8,8.8) coordinate (v2);
\path (5.2,8.8) coordinate (v3);
\path (6.8,9.5) coordinate (v4);
\path (5.2,9.5) coordinate (v5);
\path (6.8,10.5) coordinate (v6);
\path (5.2,10.5) coordinate (v7);
\path (6,12.3) coordinate (v8);

%%\path (5.2,9.5) coordinate (v9);
\path (7.5,8.5) coordinate (v9);
\path (8,7.2) coordinate (v10);
\path (8.5,9) coordinate (v11);
\path (9.1,7.7) coordinate (v12);
\path (9.7,9.5) coordinate (v13);
\path (10.3,8.2) coordinate (v14);
\path (11.5,9.5) coordinate (v15);
%%\path (4.5,8.5) coordinate (v7);
\path (2,6) coordinate (v16);
\path (2.5,4.7) coordinate (v17);
\path (3,6.5) coordinate (v18);
\path (3.6,5.2) coordinate (v19);
\path (4.2,7) coordinate (v20);
\path (4.8,5.7) coordinate (v21);
\path (0.8,4.7) coordinate (v22);
%% edges %%
\draw (v1) -- (v2) -- (v4);
\draw (v6) -- (v8) -- (v7);
\draw (v5) -- (v3) -- (v1);
\draw (v1) -- (v10) -- (v12); 
\draw (v14) -- (v15) -- (v13);
\draw (v11) -- (v9) -- (v1);
\draw  [dotted](v4) -- (v6);
\draw  [dotted](v7) -- (v5);
\draw  [dotted](v12) -- (v14);
\draw  [dotted](v13) -- (v11);
\draw (v18) -- (v20) -- (v1) -- (v21) -- (v19); 
\draw (v16) -- (v22) -- (v17);
\draw (v11) -- (v9) -- (v1);
\draw  [dotted](v16) -- (v18);
\draw  [dotted](v17) -- (v19);
\draw  [dotted] (4.6,7.5) -- (5.2,8.1);
\draw  [dotted] (5.6,6.0) -- (7.3,6.9);

%% vertices %%%
\draw (v1)  [fill=white] circle (\vr);
\draw (v2)  [fill=white] circle (\vr);
\draw (v3)  [fill=white] circle (\vr);
\draw (v4)  [fill=white] circle (\vr);
\draw (v5)  [fill=white] circle (\vr);
\draw (v6)  [fill=white] circle (\vr);
\draw (v7)  [fill=white] circle (\vr);
\draw (v8)  [fill=white] circle (\vr);
\draw (v9)  [fill=white] circle (\vr);
\draw (v10)  [fill=white] circle (\vr);
\draw (v11)  [fill=white] circle (\vr);
\draw (v12)  [fill=white] circle (\vr);
\draw (v13)  [fill=white] circle (\vr);
\draw (v14)  [fill=white] circle (\vr);
\draw (v15)  [fill=white] circle (\vr);

%%\draw (v16)  [fill=white] circle (\vr);

\draw (v16)  [fill=white] circle (\vr);
\draw (v17)  [fill=white] circle (\vr);
\draw (v18)  [fill=white] circle (\vr);
\draw (v19)  [fill=white] circle (\vr);
\draw (v20)  [fill=white] circle (\vr);
\draw (v21)  [fill=white] circle (\vr);
\draw (v22)  [fill=white] circle (\vr);
%% text %%
\draw[below] (v1) node {$u$}; 
\end{scope}

\end{tikzpicture}
\end{center}
\caption{The graph $rC_{2q}(u)$}
\label{fig:glued cycles}
\end{figure}
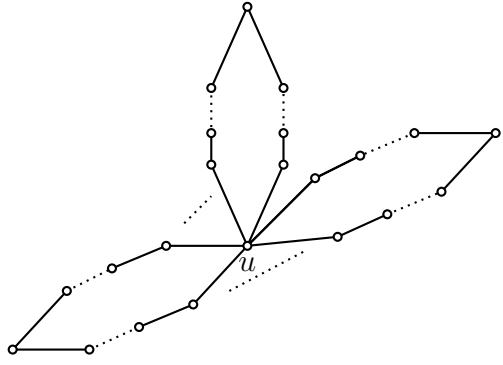

For a graph $G$ and a vertex $u\in V(G)$, let $rG(u)$ denote the graph obtained from $r$ disjoint copies of $G$ by identifying a copy of $u$ in each of the copies. In Fig.~\ref{fig:glued cycles} the graph $rC_{2q}(u)$ is schematically presented, where $u$ is an arbitrary, fixed vertex of the even cycle $C_{2q}$. 

%%%%%%%%%%%%%%%%%%%%%%%%%%%%%%%%%%%%%%%%%%%%%%%%%%

A graph $G$ is \textit{transmission regular} if all its vertices have the same transmission, cf.~\cite{abiad-2017, liu-2022, LDW2016}. Now we have: 

\begin{proposition}
\label{prop:not2-con}
Let $G$ be a bipartite, transmision regular graph, $u \in V(G)$, and $r\ge 2$. Then the graph $rG(u)$ is a $((r-1)(n-1))$-STI graph. 
\end{proposition}

\proof
let $e=vw$ be an arbitrary edge of $G$, and consider it in a copy $G'$ of $G$ in $rG(u)$. As $G$ is bipartite, we may without loss of generality assume that $d_{G'}(v,u) < d_{G'}(w,u)$.  Moreover $G$ is transmission regular it  follows that $n_v(e|G') = n_w(e|G')$. Hence 
$$I_{rG(u)}(vw)= |n_v(e|G') - n_w(e|G')| + (n-1)(r-1) = (n-1)(r-1)$$
and we are done. 
\qed

\begin{figure}[ht!]
\begin{center}
\begin{tikzpicture}[scale=0.5,style=thick]
\tikzstyle{every node}=[draw=none,fill=none]
\def\vr{2.5pt} 

\begin{scope}[yshift = 0cm, xshift = 0cm]
\path (7.5,9.5) coordinate (v0);
\path (9.5,7.5) coordinate (v2);
\path (9.5,4.5) coordinate (v3);
\path (7.5,2.5) coordinate (v4);
\path (4.5,2.5) coordinate (v5);
%%v6&\ddots&v5
\path (2.5,4.5) coordinate (v6);
\path (2.5,7.5) coordinate (v7);
\path (4.5,9.5) coordinate (v1);
\path (9.5,9.5) coordinate (u1);
\path (11,6) coordinate (u2);
\path (6,1) coordinate (u3);
\path (1,6) coordinate (u4);
\path (2.5,9.5) coordinate (u5);
\path (6,11) coordinate (w{1,1});
\path (7.5,7.5) coordinate (w1);
\path (8,6) coordinate (w2);
\path (6,4) coordinate (w3);
\path (4,6) coordinate (w4);
\path (4.5,7.5) coordinate (w5);
\path (6,8) coordinate (w{1,p});

\path (9,9) coordinate (z1);
\path (10.2,6) coordinate (z2);
\path (6,1.8) coordinate (z3);
\path (1.8,6) coordinate (z4);
\path (3,9) coordinate (z5);
\path (6,10.2) coordinate (w{1,2});
%% edges %%
\draw (v6) -- (z4) -- (v7) -- (z5) -- (v1) -- (w{1,2}) -- (v0) -- (u1) -- (v2) -- (u2) -- (v3) -- (w2) -- (v2) -- (w1) -- (v0) --(w{1,p}) -- (v1) -- (w5) -- (v7) -- (w4) -- (v6) -- (u4) -- (v7) -- (u5) -- (v1) -- (w{1,1}) -- (v0) -- (z1) -- (v2) --(z2) -- (v3);  
\draw (v4) -- (u3) -- (v5) -- (w3) -- (v4) -- (z3) -- (v5);
%%\draw  (8.5, 3.5) circle (1cm);
\draw [dotted] (8.6,8.6) -- (7.9,7.9);
\draw [dotted] (9.6,6) -- (8.6,6);
\draw [dotted] (6,2.4) -- (6,3.4);
\draw [dotted] (2.4,6) -- (3.4,6);
\draw [dotted] (3.4,8.6) -- (4.1,7.9);
\draw [dotted] (6,9.6) -- (6,8.6);
\draw [dotted] (8.9,3.9) -- (8.1,3.1);
\draw [dotted] (3.9,3.1) -- (3.1,3.9);
%% vertices %%%

%%\draw (1,2) .. controls (.2.5,3) .. (1,4);
\draw (v0)  [fill=white] circle (\vr);
\draw (v2)  [fill=white] circle (\vr);
\draw (v3)  [fill=white] circle (\vr);
\draw (v4)  [fill=white] circle (\vr);
\draw (v5)  [fill=white] circle (\vr);
\draw (v6)  [fill=white] circle (\vr);
\draw (v7)  [fill=white] circle (\vr);
\draw (v1)  [fill=white] circle (\vr);
\draw (u1)  [fill=white] circle (\vr);
\draw (u2)  [fill=white] circle (\vr);
\draw (u3)  [fill=white] circle (\vr);
\draw (u4)  [fill=white] circle (\vr);
\draw (u5)  [fill=white] circle (\vr);
\draw (w{1,1})  [fill=white] circle (\vr);
\draw (w1)  [fill=white] circle (\vr);
\draw (w2)  [fill=white] circle (\vr);
\draw (w3)  [fill=white] circle (\vr);
\draw (w4)  [fill=white] circle (\vr);
\draw (w5)  [fill=white] circle (\vr);
\draw (w{1,p})  [fill=white] circle (\vr);

\draw (z1)  [fill=white] circle (\vr);
\draw (z2)  [fill=white] circle (\vr);
\draw (z3)  [fill=white] circle (\vr);
\draw (z4)  [fill=white] circle (\vr);
\draw (z5)  [fill=white] circle (\vr);
\draw (w{1,2})  [fill=white] circle (\vr);

%% text %%
\draw[above] (w{1,1}) node {$w_{1,1}$};
\draw[below] (w{1,2}) node {$w_{1,2}$};
\draw[below] (w{1,p}) node {$w_{1,p}$};
\draw[above] (v1) node {$v_{1}$};
\end{scope}
\end{tikzpicture}
\end{center}
\caption{The graph $\Gamma_{p,q}$ }
\label{fig:Gpq}
\end{figure}
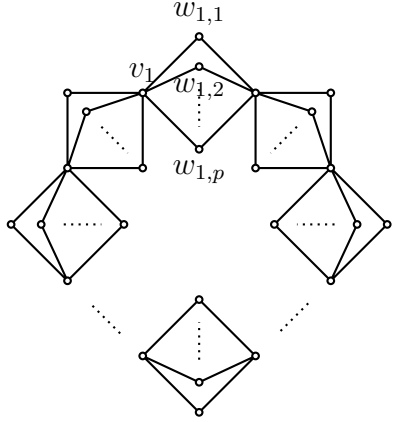

If $p,q \ge 2$, then let $\Gamma_{p,q}$ be the graph with the vertex set $V(\Gamma_{p,q}) = \{v_i , w_{i,j}:\ i\in [2q], j \in [p]\}$ and the edge set 
\begin{align*}
E(\Gamma_{p,q}) = &  \{v_iw_{i,j}, v_iw_{i-1,j}:\ 2 \le i \le 2q-1 , j\in [p]\}
\cup \{ v_iw_{i,j}:\  i\in \{1,2q\},  j\in [p]\}\\ 
& \cup\{v_1w_{2q,j}, v_{2q}w_{2q-1,j}:\ j\in [p] \}\,.
\end{align*}
To put it more informally, $\Gamma_{p,q}$ is obtained from $2q$ copies of $K_{2,p}$ by attaching them in circular manner, see Fig.~\ref{fig:Gpq}. 

%%%%%%%%%%%%%%%%%%%%%%%%%%%%%%%%%%%%%%%%%

\begin{proposition}
If $p,q \ge 2$,  then $\Gamma_{p,q} $ is  a $(2p-2)$-STI graph.
\end{proposition}

\proof
By the symmetry of the graph $\Gamma_{p,q}$ we infer that $ I(v_iw_{i,j})$ is independent of the section of $i$ and $j$. Hence it suffices to compute $I(v_1w_{1,1})$. Setting $e = v_1w_{1,1}$ we have  
\begin{align*}
N_{v_1}(e) = \{v_1\} & \cup \{ w_{1,j}:\  2\le j \le p \} \cup \{v_i:\ q+2 \le i \le 2q \} \\
& \cup \{w_{i,j}:\  q+1 \le i \le 2q, j \in [p] \}\,.
\end{align*}
Therefore, $n_{v_1}(e) = 1 + (p-1) + (q-1) + (pq) = pq + p + q-1$ and $n_{w_{1,1}}(e)= 2q(p+1) - (pq+p+q-1)= q(p+1) - p+1$. By Theorem~\ref{thm:Entringer} we conclude that  
\begin{align*}
 |I(v_1w_{1,1})| & = |n_{v_1} - n_{w_{1,1}}| \\ 
  & = |(pq+p+q-1) - (pq + q -p+1)|\\
  & = 2p-2\,, 
\end{align*}
hence the assertion. 
\qed

\begin{figure}[ht!]
\begin{center}
\begin{tikzpicture}[scale=0.7,style=thick]
\tikzstyle{every node}=[draw=none,fill=none]
\def\vr{2.5pt} 

\begin{scope}[yshift = 0cm, xshift = 0cm]
\path (7.5,9.5) coordinate (v1);
\path (9.5,8) coordinate (v2);
\path (11,6) coordinate (v3);
\path (9.5,4.8) coordinate (v4);
\path (7.5,3.3) coordinate (v5);
\path (5.5,4.8) coordinate (v6);
\path (4,6) coordinate (v7);
\path (5.5,8) coordinate (v8);
\path (9.3,9.5) coordinate (u1);
\path (11,7.8) coordinate (u2);
\path (9.3,3.3) coordinate (u3);
\path (5.7,3.3) coordinate (u4);
\path (4,7.8) coordinate (u5);
\path (5.7,9.5) coordinate (u6);
\path (7.9,8) coordinate (w1);
\path (9.5,6.4) coordinate (w2);
\path (7.9,4.8) coordinate (w3);
\path (7.1,4.8) coordinate (w4);
\path (5.3,6.3) coordinate (w5);
\path (7.1,8) coordinate (w6);

\path (8.9,9.1) coordinate (z1);
\path (8.9,3.7) coordinate (z2);
\path (4.4,7.3) coordinate (z3);

%% edges %%
\draw (v1) -- (u1) -- (v2) -- (u2) -- (v3) -- (w2) -- (v2) -- (w1) -- (v1) -- (w6) -- (v8) -- (w5) -- (v7) -- (u5) -- (v8) --(u6) -- (v1); %%   
\draw (v4) -- (u3) -- (v5) -- (u4) -- (v6) -- (w4) -- (v5) -- (w3) -- (v4);
\draw (v1) -- (z1) -- (v2);
\draw (v4) -- (z2) -- (v5);
\draw (v7) -- (z3) -- (v8);

\draw  [dotted](10.6,5.5) -- (9.9,5.1);
\draw  [dotted](5.1,5.1) -- (4.4,5.7);
\draw  [dotted](8.6,8.8) -- (8.2,8.3);
\draw  [dotted](8.6,4.0) -- (8.1,4.5);
\draw  [dotted](4.6,7.05) -- (5.1,6.55);

%% vertices %%%

\draw (v1)  [fill=white] circle (\vr);
\draw (v2)  [fill=white] circle (\vr);
\draw (v3)  [fill=white] circle (\vr);
\draw (v4)  [fill=white] circle (\vr);
\draw (v5)  [fill=white] circle (\vr);
\draw (v6)  [fill=white] circle (\vr);
\draw (v7)  [fill=white] circle (\vr);
\draw (v8)  [fill=white] circle (\vr);
\draw (u1)  [fill=white] circle (\vr);
\draw (u2)  [fill=white] circle (\vr);
\draw (u3)  [fill=white] circle (\vr);
\draw (u4)  [fill=white] circle (\vr);
\draw (u5)  [fill=white] circle (\vr);
\draw (u6)  [fill=white] circle (\vr);
\draw (w1)  [fill=white] circle (\vr);
\draw (w2)  [fill=white] circle (\vr);
\draw (w3)  [fill=white] circle (\vr);
\draw (w4)  [fill=white] circle (\vr);
\draw (w5)  [fill=white] circle (\vr);
\draw (w6)  [fill=white] circle (\vr);

\draw (z1)  [fill=white] circle (\vr);
\draw (z2)  [fill=white] circle (\vr);
\draw (z3)  [fill=white] circle (\vr);

%% text %%
\draw[above] (v1) node {$v_{1}$};
\draw[ right] (v2)++(0,0.3) node {$v_{2}$};
\draw[left] (v8)++(0,0.3) node {$v_{2q}$};
\draw[above] (u1) node {$w_{1,1}$};
\draw[above] (u6) node {$w_{2q,1}$};
\draw[below] (w1)++(0.2,0) node {$w_{1,p}$};
\draw[below] (w6)++(-0.3,0) node {$w_{2q,2}$};
\draw[above] (u2) node {$w_{2,1}$};
\draw[left] (w2) node {$w_{2,2}$};

\end{scope}

\end{tikzpicture}
\end{center}
\caption{The graph $H_{p,q}$ }
\label{fig3}
\end{figure}
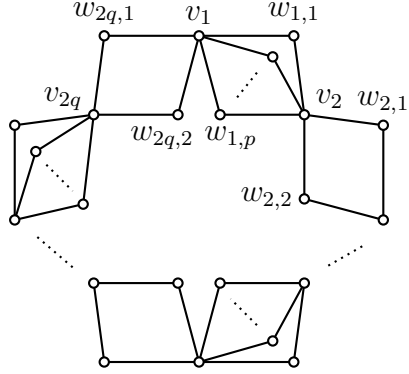

If $p,q \ge 2$, then let $H_{p,q}$  be the graph with the vertex set 
$$V(H_{p,q}) = \{v_i:\   i \in [2q] \} \cup \{w_{2r-1,j}:\ r \in [q], j \in [p] \} \cup \{w_{2r,j}:\ r \in [q], j \in [2] \}$$
and the edge set $A \cup B \cup C$, where
\begin{align*}
A & =  \{v_1w_{2q,j}:\ j \in [2] \}\,,\\
B & = \{v_{2r-1}w_{2r-1,j}:\ r\in [q], j\in [p] \} \cup \{v_{2r-1}w_{2r-2,j}:\ 2\le r\le q, j \in [2] \}\,,\\
C & = \{v_{2r}w_{2r,j}:\  r \in [q], j\in [2] \} \cup \{v_{2r}w_{2r-1,j}:\ r\in [q], j\in [p] \}\,. 
\end{align*}
Informatively, $H_{p,q}$ has a cyclic structure where $K_{2,p}$ and $K_{2,2}$ alternate, see Fig.~\ref{fig3}.

\begin{proposition}
If $p \ge 2$ and $q\ge 3$ is odd, then $H_{p,q}$ is a $p$-STI graph.
\end{proposition}

\proof
For the rest of the proof set $n = n(H_{p,q}) = q(p+4)$, Considering the symmetry of $H_{p,q}$, it is enough to determine the transmission imbalance for the edges $e_1= v_1w_{1,1}$ and $f_1=v_1w_{2q,1}$. For the edge  $e_1=v_1w_{1,1}$ we have
\begin{align*}
N_{w_{1,1}}(e_1) = & \{w_{1,1} \} \cup \{v_i:\  2\le i \le q+1\} \cup \{w_{r,j}:\ r\ \mbox{even}, 2\le r \le q, j \in 2 \}\\
& \cup \{w_{r,j}:\ r\ \mbox{odd}, 3\le r \le q, j \in p \}\,.
\end{align*}
Therefore $n_{w_{1,1}}(e_1) = (q+1)+2(\frac{q-1}{2}) + p(\frac{q-1}{2})$ and hence
$$I(e_1)= |n - 2n_{w_{1,1}}(e)| = |q(p+4) -2(q+1+(p+2)(\frac{q-1}{2}))| = p\,.$$
For the edge $ f_1 = v_{1}w_{2q,1}$ we have
\begin{align*}
N_{v_1}(f_1)  = & \{w_{2q,2} \} \cup \{v_i:\ i \in [q] \} \cup \{w_{r,j}:\ r\ \mbox{even}, 2\le r \le q-1, j \in [2] \} \\
 & \cup \{w_{r,j}:\ r\ \mbox{odd}, r \in [q], j \in [p] \}\,.
\end{align*}
Thus $n_{v_1}(f_1)= 1 + q + 2\frac{q-1}{2} + p\frac{q+1}{2}= 2q-1+p\frac{q+1}{2}$ and we obtain that 
$$I(f_1)= |n-2n_{v_1}(f_1)| =|q(p+4) - 2(2q-1+p\frac{q+1}{2})| = p\,.$$
We conclude that $H_{p,q}$ is a $p$-STI graph. 
\qed

We add that if $q$ is even, then the transmission imbalances of the edges of $H_{p,q}$ are $2p-2$ and $2$. 

\begin{figure}[ht!]
\begin{center}
\begin{tikzpicture}[scale=0.8,style=thick]
\tikzstyle{every node}=[draw=none,fill=none]
\def\vr{2.5pt} 

\begin{scope}[yshift = 0cm, xshift = 0cm]
\path (1,4) coordinate (v{n-r});
\path (2.5,4.5) coordinate (v{n-r+1});
\path (4,4) coordinate (v1);
\path (5,3) coordinate (v4);
\path (5,1.5) coordinate (v5);
\path (4,0.5) coordinate (v6);
\path (2.5,0) coordinate (v7);
\path (1,0.5) coordinate (v8);
\path (0,1.5) coordinate (v9);
\path (0,3) coordinate (v10);
\path (2.5,3.7) coordinate (v{n-1});
\path (2.5,3.4) coordinate (vn);

%% edges %%
\draw (v{n-r}) -- (v{n-r+1}) -- (v1) -- (v4) -- (v5) -- (v6);
\draw [dotted](v6) -- (v7);
\draw [dotted](v7)-- (v8);
\draw (v8)-- (v9) -- (v10) -- (v{n-r}) -- (vn) -- (v1);
\draw [dotted](v{n-r+1})-- (v{n-1});
\draw (v{n-r}) -- (v{n-1}) -- (v1);
%% vertices %%%
\draw (v{n-r})  [fill=white] circle (\vr);
\draw (v{n-r+1})  [fill=white] circle (\vr);
\draw (v1)  [fill=white] circle (\vr);
\draw (v4)  [fill=white] circle (\vr);
\draw (v5)  [fill=white] circle (\vr);
\draw (v6)  [fill=white] circle (\vr);
\draw (v7)  [fill=white] circle (\vr);
\draw (v8)  [fill=white] circle (\vr);
\draw (v9)  [fill=white] circle (\vr);
\draw (v10)  [fill=white] circle (\vr);
\draw (v{n-1})  [fill=white] circle (\vr);
\draw (vn)  [fill=white] circle (\vr);
%% text %%
\draw[left] (v{n-r}) node {$v_{n-r}$};

\draw[above] (v{n-r+1}) node {$v_{n-r+1}$};
\draw[right] (v1) node {$v_1$};
\draw[below] (vn) node {$v_{n}$};
\end{scope}

\end{tikzpicture}
\end{center}
\caption{The graph $G_{n,r}$}
\label{fig:Gnr}
\end{figure}
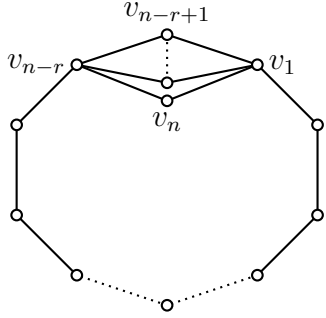

The last class of graphs we present is defined as follows. If $r\ge 2$ and $n \ge 5$, then let $G_{n,r}$ be the graph with the vertex set $\{v_1, v_2, \ldots, v_n \}$ and the edge set $\{ v_iv_{i+1}:\  i\in [n-r-1]\} \cup \{v_1v_j, v_{n-r}v_j:\ n-r+1 \le j \le n\}$, see Fig.~\ref{fig:Gnr}. 

\begin{proposition}
If $n\ge 5$, $r\ge 2$, and $n \equiv r-1 \pmod {2}$, then $G_{n,r}$ is a $(r-1)$-STI graph. 
\end{proposition}
\proof
Let $v_iv_j$ be an edge of $G = G_{n,r}$, where $v_i$ is closer to one of two vertices $v_1$ and $v_{n-r}$ than to $v_j$, that is, $\min\{d_G(v_i,v_1), d_G(v_i,v_{n-r})\} < \min\{d_G(v_j,v_1), d_G(v_j,v_{n-r})\}$. Note that the cycle $C: v_1 v_2 \cdots v_{n-r}v_tv_1$, where $n-r+1 \le t \le v_n$, is an even cycle and then $n(v_i|C) = n(v_j|C)$. Moreover the other vertices not on the cycle $C$ are closer to $v_i$ than to $v_j$. 
 Thus
$$ I(v_iv_j) =|n_{v_i}-n_{v_j}|= \left|\left(\frac{n-r+1}{2} + r-1 \right) -\frac{n-r+1}{2}\right| = r-1$$
and we are done. 
\qed

%%%%%%%%%%%%%%%%%%%%%%%%%%%%%%%%%%%%%%%%%%%%%%%%%%%%%%%
\section*{Acknowledgements}
%%%%%%%%%%%%%%%%%%%%%%%%%%%%%%%%%%%%%%%%%%%%%%%%%%%%%%%

Sandi Klav\v{z}ar was supported by the 
% Javna agencija za znanstvenoraziskovalno in inovacijsko dejavnost Republike Slovenije (ARIS)
Slovenian Research Agency (ARIS) under the grants P1-0297, J1-2452, and N1-0285.

\end{document}